\documentclass[11pt]{article}
\usepackage{numbysec,amsfonts,amssymb,amsmath,amsthm}

-\marginparwidth \oddsidemargin -4mm \evensidemargin \oddsidemargin

\usepackage{graphicx,epsf}

\advance\hoffset by 5mm
\newtheorem{thm}[theorem]{Theorem}
\newtheorem{prop}[theorem]{Proposition}
\newtheorem{lemma}[theorem]{Lemma}

\newtheorem{corollary}[theorem]{Corollary}
\theoremstyle{definition} \newtheorem{example}[theorem]{Example}
\newcommand{\Rm}{{\mathbb R}}
\newcommand{\pdr}[2]{\frac{\partial{#1}}{\partial{#2}}}

\newcommand{\eps}{\varepsilon}

\newcommand{\Tm}{{\mathbb T}}
\newcommand{\Zm}{{\mathbb Z}}
\newcommand{\calS}{{\mathcal S}}

\newcommand{\vx}{{\bf x}}

\newcommand{\calC}{{\mathcal{C}}}
\newcommand{\commentout}[1]{}

\newcommand{\cea}{c_e^*(A)}
\newcommand{\dea}{D_e(A)}
\numberbysection

\begin{document}
\title{KPP Pulsating Front Speed-up by Flows}
\author{
Lenya Ryzhik\thanks{Department of Mathematics, University of
Chicago, Chicago, IL 60637, USA; ryzhik@math.uchicago.edu} \and
Andrej Zlato\v s\thanks{Department of Mathematics, University of
Chicago, Chicago, IL 60637, USA; zlatos@math.uchicago.edu} }

\maketitle

\begin{abstract}
We obtain a criterion for pulsating front speed-up by general
periodic incompressible flows in two dimensions and in the presence
of KPP nonlinearities. We achieve this by showing that the ratio of
the minimal front speed and the effective diffusivity of the flow is
bounded away from zero and infinity by constants independent of the
flow. We also study speed-up of reaction-diffusion fronts by various
examples of flows in two and three dimensions.
\end{abstract}

\section{Introduction}

We consider reaction-diffusion fronts propagating in a strong
periodic incompressible flow on $\Rm^n$:
\begin{equation}\label{intro-pf-1}
T_t+Au\cdot\nabla T=\Delta T+f(T).
\end{equation}
Here $T(t,x)\in[0,1]$ is the normalized temperature and the
nonlinearity $f$ is of the KPP type: $f(s)$ is a Lipschitz function
such that $f(0)=f(1)=0$, while $f(s)>0$  and $f(s)\le f'(0)s$ for
$s\in(0,1)$. The Lipschitz flow $u(x)$ is $1$-periodic,
incompressible, and has mean zero. That is, $u(x+k)=u(x)$ when
$k\in\Zm^n$, $\nabla\cdot\ u=0$, and
\[
\int_{{\mathbb T}^n}u(x)dx=0.
\]
The parameter $A\in\Rm$ is the amplitude of the advection, and we
will mainly be interested in $A\gg 1$.

It has been proved in \cite{BH} that when $u(x)$ is periodic,
equation (\ref{intro-pf-1}) has pulsating front solutions of the
form $T(t,x)=U(x\cdot e-ct,x)$, where $c>0$ is the propagation speed
and $e\in\Rm^n$ is the unit vector in the direction of propagation.
The function $U(s,x)$ is periodic in $x\in\Rm^n$ and has uniform in
$x\in\Tm^n$ limits as $s\to\pm\infty$:
\begin{eqnarray}\label{intro-pf-2}
&&\lim_{s\to-\infty}U(s,x)=1,\\
&&\lim_{s\to+\infty}U(s,x)=0. \nonumber
\end{eqnarray}
Pulsating front solutions were shown in \cite{BH} to exist for all
$|e|=1$ and all $c\ge \cea$. As in the one-dimensional case without
advection \cite{KPP}, the minimal front speed $\cea$ (we suppress
the $u$ and $f$ dependence in our notation) determines the
propagation speed of solutions of the Cauchy problem for
(\ref{intro-pf-1}) with general compactly supported initial data,
and is therefore of a special interest \cite{BHN-1,Weinberger}.

The presence of an incompressible flow in (\ref{intro-pf-1})
improves mixing due to diffusion and is thus expected to enhance the
speed of reaction-diffusion fronts. This problem has been studied
actively in the recent years, especially in the large $A$ limit: it
has been shown in \cite{B,CKOR,Heinze} that the pulsating front
speed in the direction of a mean-zero shear flow behaves as
$\cea=O(A)$ for large $A$, and in \cite{NR} that $\cea=O(A^{1/4})$
for cellular flows in two dimensions. In both of these cases the
minimal front speed scales as $\cea\sim\sqrt{\dea}$ for $A\gg 1$
where $\dea$ is the corresponding effective diffusivity of the flow
$Au$ in the direction $e$. In the present paper we use the method of
\cite{NR} to show that this is indeed the case in general in two
dimensions.



\commentout{


 Let $u$ be an incompressible 1-periodic mean-zero
flow, let $f$ be KPP, and $e\in\Rm^2$ a unit vector. We consider
the pulsating fronts in the direction $e$, that is solutions of
\begin{equation}\label{pf-1}
T_t+Au\cdot\nabla T=\Delta T+f(T)
\end{equation}
of the form $T(t,x_1,x_2)=U(x\cdot e-ct,x_1,x_2)$, where $c>0$ and
\begin{eqnarray}\label{pf-2}
&&U(s,x_1,x_2)=U(s,x_1,x_2+1),\nonumber\\
&&U(s,x_1,x_2)=U(s,x_1+1,x_2),\nonumber\\
&&\lim_{s\to-\infty}U(s,x_1,x_2)=1,\nonumber\\
&&\lim_{s\to+\infty}U(s,x_1,x_2)=0. \nonumber
\end{eqnarray}
The last two limits are required to be uniform in $s$.   It has been
shown in \cite{BH} that pulsating front solutions exist for all
$c\ge \cea$. We are interested in the asymptotics of the minimal
front speed $\cea$ as $A\to+\infty$.

}


Let us recall the definition of the effective diffusivity. Consider
the advection-diffusion problem
\begin{equation}\label{ed-1}
p_t+Au\cdot\nabla p=\Delta p
\end{equation}
with $u$ a periodic incompressible flow. The long-time behavior of
the solutions of \eqref{ed-1} is governed by the effective diffusion
equation
\begin{equation}\label{ed-12}
\bar p_t=\sum_{i,j=1}^n \sigma_{ij}(A)\frac{\partial^2\bar
p}{\partial x_i\partial x_j}.
\end{equation}
The ($x$-independent) effective diffusivity  matrix $\sigma(A)$ is
obtained as follows. For any $e\in\Rm^n$, let $\chi_e(x)$ be the
periodic mean-zero solution of the cell problem
\begin{equation}\label{ed-3}
-\Delta\chi_e+Au\cdot\nabla\chi_e=Au\cdot e
\end{equation}
on $\Tm^n$. Then the matrix $\sigma(A)$ is given by
\begin{equation}\label{ed-4}
e\cdot\sigma(A)e' = \int_{\Tm^n}(\nabla\chi_e +
e)\cdot(\nabla\chi_{e'} + e') dx = e\cdot e' +
\int_{\Tm^n}\nabla\chi_e\cdot\nabla\chi_{e'} dx,
\end{equation}
for any $e,e'\in\Rm^n$. The effective spreading in the direction $e$
is then governed by the effective diffusivity
\[
\dea=e\cdot\sigma(A)e = 1+\int_{\Tm^n}|\nabla\chi_e|^2 dx.
\]

When the nonlinearity in \eqref{intro-pf-1} is weak and
\eqref{intro-pf-1} becomes
\[
T_t+Au\cdot\nabla T=\Delta T+\eps f(T)
\]
with $\eps\ll 1$, one may consider the long time--large space
scaling $t\to t/\eps^2$, $x\to x/\eps$ leading to
\[
T_t+\frac{A}{\eps}u\left(\frac{x}{\eps}\right)\cdot\nabla T=\Delta
T+f(T).
\]
The homogenized version of this equation is
\[
T_t=\nabla\cdot(\sigma(A)\nabla T)+f(T),
\]
with the corresponding homogenized minimal front speed
$v_e^*(A)=2\sqrt{f'(0)\dea}$. This approximation  holds only on
certain time--space scales in the original variables (namely,
$t=O(1/\eps^2)$ and $x=O(1/\eps)$). However, it suggests a relation
between the minimal front speed for \eqref{intro-pf-1} and the
effective diffusivity. The following result confirms this relation
in two dimensions.


\begin{thm}\label{thm-main}
There exists $C>0$ (independent of $A,u,f,e$) such that if $u(x)$ is
a 1-periodic incompressible Lipschitz flow on $\Rm^2$, $f$ a KPP
nonlinearity, $A\in\Rm$, and $e\in\Rm^2$ a unit vector, then
\begin{equation}\label{intro-D-c}
\frac{\sqrt{f'(0)}}{C\big(1+\sqrt{f'(0)}\,\big)}\le
\frac{\cea}{\sqrt{\dea}}\le C\sqrt{f'(0)}\big(1+\sqrt{f'(0)}\,\big).
\end{equation}
Moreover, there is $f_0>0$ such that
\begin{equation}\label{eq-thm-ratio}
\bigg|\frac{\cea}{2\sqrt{f'(0)\dea}} - 1 \bigg| \le C\,f'(0)^{1/4}.
\end{equation}
whenever $0<f'(0)\le f_0$.
\end{thm}

That is, the ratio $\cea/\sqrt{\dea}$ is bounded away from zero and
infinity by constants only dependent on $f'(0)$, and becomes close
to $2\sqrt{f'(0)}$ when $f'(0)$ is small. We note that the slightly
weaker, and only upper bound $\cea/\sqrt{\dea} \le C_\eps
\sqrt{f'(0)}(1+\sqrt{f'(0)})\|Au\|_\infty^{\eps}$ for any $\eps>0$
has been obtained in \cite{Heinze}.

As $f'(0)\to+\infty$, the lower bound in \eqref{intro-D-c} stays
bounded whereas the upper one grows linearly with $f'(0)$. We show
by looking at the example of shear flows (see Example \ref{E.4.1a})
below that at least the lower bound cannot be improved. This conclusion
can also be reached using the results of \cite{Heinze} for the shear flow.

It follows from Theorem \ref{thm-main} that the minimal front speed
$\cea$ has the same asymptotic behavior in the regime of large $A$
as does $\sqrt{\dea}$. But for the latter quantity we have the
following general sharp criterion which holds in any spatial
dimension.

\begin{prop}\label{thm-main0}
Let $u(x)$ be a 1-periodic incompressible Lipschitz flow on $\Rm^n$
and let $e\in\Rm^n$ be a unit vector.

(i) If the equation
\begin{equation}\label{eq-thm-eigen-fn}
u\cdot\nabla\phi_e=u\cdot e
\end{equation}
has a solution $\phi_e\in H^1(\Tm^n)$, then
\begin{equation}\label{intro-bdd}
\limsup_{A\to+\infty}\dea<+\infty.
\end{equation}

(ii) If (\ref{eq-thm-eigen-fn}) has no $H^1(\Tm^n)$-solutions, then
\begin{equation}\label{intro-infty}
\lim_{A\to+\infty}\dea=+\infty.
\end{equation}
\end{prop}

{\it Remark.} Note that it follows from part (i) that the set of all non-negative
multiples of unit vectors $e\in\Rm^n$ for which \eqref{intro-bdd}
holds is a subspace of $\Rm^n$. Indeed, the sum of two solutions of
\eqref{eq-thm-eigen-fn} for $e$ and $e'$ is a solution for $e+e'$,
and the negative of a solution for $e$ is a solution for $-e$.
\medskip

The result of Proposition~\ref{thm-main0} is not new. It has already
appeared in \cite{FP}, although it has been stated only in two
dimensions, and the first claim has also appeared earlier in
\cite{BGW}. Much more precise asymptotic behavior of $\dea$ is well
understood for many specific examples of flows --- see \cite{KM} for
an extensive list of references.

Putting Theorem \ref{thm-main} and Proposition \ref{thm-main0}
together, we have the following characterization of flows in two
dimensions that speed up KPP fronts.

\begin{corollary} \label{cor-main}
Let $u(x)$ be a 1-periodic incompressible Lipschitz flow on $\Rm^2$,
let $f$ be a KPP nonlinearity, and $e\in\Rm^2$ a unit vector.

(i) If \eqref{eq-thm-eigen-fn} has a solution $\phi_e\in
H^1(\Tm^2)$, then
\begin{equation}\label{intro-bdd2}
\limsup_{A\to+\infty}\cea<+\infty.
\end{equation}

(ii) If \eqref{eq-thm-eigen-fn} has no $H^1(\Tm^2)$-solutions, then
\begin{equation}\label{intro-infty2}
\lim_{A\to+\infty}\cea=+\infty.
\end{equation}
\end{corollary}


{\it Remark.} In particular, the pulsating front speed for KPP
nonlinearities may not diverge to $+\infty$ along some sequence of
amplitudes while staying bounded along another sequence.
\medskip


The paper is organized as follows. Section~\ref{sec:proof} contains
the proofs of Theorem~\ref{thm-main} and
Proposition~\ref{thm-main0}. The proof of Theorem~\ref{thm-main} is
based on the proof of the main result of \cite{NR}.
Section~\ref{sec:gen} contains the generalization of Corollary
\ref{cor-main}(ii) to higher dimensions and various examples. In
particular, we show there that the minimal front speed for a class
of cellular three-dimensional flows satisfies \eqref{intro-infty2}.
To the best of our knowledge, this is the first time that the
front speed-up
by a cellular three-dimensional flow has been established.

{\bf Acknowledgment.} This work has been supported by ASC Flash
Center at the University of Chicago.  LR was supported by NSF
grant DMS-0604687 and AZ by NSF grant DMS-0632442.

\section{Diffusivity enhancement and front speed-up}\label{sec:proof}


We first present the proof of Proposition~\ref{thm-main0} for the
convenience of the reader.

\commentout{ 

\begin{prop}\label{prop-speed-finite}
Assume that (\ref{eq-thm-eigen-fn}) has no $H^1(\Tm^2)$
solution. Then the effective diffusivity satisfies
\begin{equation}\label{eq-prop-infinite}
\lim_{A\to+\infty}\dea=+\infty.
\end{equation}
\end{prop}

\begin{prop}\label{prop-eff-diff-finite}
Assume that (\ref{eq-thm-eigen-fn}) has an $H^1(\Tm^2)$
solution. Then the effective diffusivity satisfies
\begin{equation}\label{eq-prop-finite}
\limsup_{A\to+\infty}\dea<+\infty.
\end{equation}
\end{prop}

\begin{prop}\label{prop-eff-diff-speed}
Let $A_n\to +\infty$, then
\begin{equation}\label{eq-prop-eff-d-s-1}
\lim_{n\to+\infty}D_e(A_n)=+\infty\qquad \hbox{if and only if}\qquad
\lim_{n\to+\infty}c^*_e(A_n)=+\infty,
\end{equation}
and
\begin{equation}\label{eq-prop-eff-d-s-2}
\limsup_{n\to+\infty}D_e(A_n)<+\infty\qquad \hbox{if and only
if}\qquad \limsup_{n\to+\infty}c^*_e(A_n)<+\infty.
\end{equation}
Moreover, (\ref{eq-thm-ratio}) holds.
\end{prop}

{\it Remark.} 1. Propositions \ref{prop-speed-finite} and
\ref{prop-eff-diff-finite} have already appeared in \cite{FP} (and
the second of them earlier in \cite{BGW}) but we prove them here for
the sake of completeness.

2. Our proofs of Propositions \ref{prop-speed-finite} and
\ref{prop-eff-diff-finite} hold in any spatial dimension (see
Theorem~\ref{T.4.1} below). However, the proof of Proposition
\ref{prop-eff-diff-speed} relies on the fact that we are working in
two dimensions.
\smallskip

}

\subsection{The proof of Proposition \ref{thm-main0}}

Let us assume that there exists a sequence $A_n\to +\infty$ and a
constant $M>0$ such that $D_e(A_n)< M$ for all $n$. It follows from
\eqref{ed-4} that there exists a sequence of mean-zero functions
$\chi_n(x)$ on $\Tm^n$ which satisfy
\begin{equation}\label{psf-1}
-\Delta\chi_n+A_nu\cdot\nabla\chi_n=A_nu\cdot e
\end{equation}
such that $\|\nabla\chi_n\|_2^2\le M$ for all $n$. As the
functions $\chi_n$ are uniformly bounded in $H^1(\Tm^n)$ there
exists a subsequence $\chi_{n_k}$ which converges to a function
$\bar\chi(x)$ weakly in $H^1(\Tm^n)$ and strongly in $L^2(\Tm^n)$
as $k\to +\infty$. We divide (\ref{psf-1}) by $A_{n_k}$ and pass
to the limit $k\to +\infty$ to obtain
\begin{equation}\label{psf-2}
u\cdot\nabla\bar\chi=u\cdot e
\end{equation}
in the sense of distributions. Since $\bar\chi\in H^1(\Tm^n)$,
(\ref{psf-2}) holds almost everywhere on $\Tm^n$. This proves
Proposition~\ref{thm-main0}(ii).

In order to prove Proposition~\ref{thm-main0}(i) let us assume that
a mean-zero function $\phi_e\in H^1(\Tm^n)$ satisfies
(\ref{eq-thm-eigen-fn}) and let $\chi_e$ be the mean-zero solution
of \eqref{ed-3}. 
Consider the function $\eta=\chi_e-\phi_e$ which satisfies
\[
-\Delta(\eta-\phi_e)+Au\cdot\nabla\eta=0.
\]
Multiplying the last equation by $\eta$ and integrating by parts,
from incompressibility of $u$ we obtain
\[
\int|\nabla\eta|^2dx=-\int\nabla\phi_e\cdot\nabla\eta dx.
\]
It follows that $\|\nabla\eta\|_2\le\|\nabla\phi_e\|_2$, and so
\[
D_e(A) = 1+ \int_{\Tm^n}|\nabla\chi_e|^2dx \le 1+
2\int_{\Tm^n}|\nabla\eta|^2dx+2\int_{\Tm^n}|\nabla\phi_e|^2dx \le 1+
4\|\nabla\phi_e\|_2^2.
\]
Therefore, $D_e(A)$ is uniformly bounded in $A$ and (i) follows.
$\Box$

\subsection{A variational principle for $c^*_e(A)$}

The proof of Theorem~\ref{thm-main} relies on a variational
principle for the effective speed which we now recall. Details and
proofs of Propositions \ref{prop-var-pr} and \ref{lem-mu-convex}
below can be found in \cite{BH,BHN-1,BHN-2,F1,NR}. Consider the
eigenvalue problem on $\Tm^n$
\begin{equation}\label{2-eig1}
\Delta\varphi-Au\cdot\nabla\varphi-2\lambda e\cdot\nabla\varphi+
{\lambda}Au\cdot e\varphi=\kappa_e(\lambda;A)\varphi,\qquad
\hbox{$\varphi>0$}.
\end{equation}
It has a unique eigenvalue $\kappa_e(\lambda;A)$ that corresponds to
a positive periodic eigenfunction $\varphi_e(x;\lambda,A)$.
\begin{prop}\label{prop-var-pr}
The minimal front
speed is described by the variational principle
\begin{equation}\label{2-c-min}
c^*_e(A)=\inf_{\lambda>0}\frac{f'(0)+\lambda^2+\kappa_e(\lambda;A)}{\lambda}.
\end{equation}
\end{prop}
It is convenient to rewrite the eigenvalue problem (\ref{2-eig1}) in
terms of the function
\[
\psi_e(x)=\varphi_e(x)e^{\lambda x\cdot e}
\]
on $\Rm^n$, with $\varphi_e$ periodically extended from $\Tm^n$.
This function is not periodic but rather belongs to the set
\[
E_{e,\lambda}^+=\left\{\hbox{$\psi(x) \,\big|\,\psi(x)e^{-\lambda
x\cdot e}$ is 1-periodic in $x$ and $\psi>0$}\right\}.
\]
The corresponding eigenvalue problem for $\psi_e(x)$ is
\begin{equation}\label{2-eig2}
{\cal L}\psi:=\Delta\psi-Au\cdot\nabla\psi=\mu_e(\lambda;A)\psi,
\qquad \hbox{$\psi\in E_{e,\lambda}^+$}
\end{equation}
with $\mu_e(\lambda;A)=\lambda^2+\kappa_e(\lambda;A)$ the unique
eigenvalue of (\ref{2-eig2}). The variational principle
(\ref{2-c-min}) may now be restated as
\begin{equation}\label{2-c-min-mu}
c^*_e(A)=\inf_{\lambda>0}\frac{f'(0)+\mu_e(\lambda;A)}{\lambda}.
\end{equation}
Let us recall some basic properties of the function
$\mu_e(\lambda;A)$. These can be found, for instance, in
\cite[Proposition 5.7]{BH} (with our $\mu_e$ being their $-h$):

\begin{prop}\label{lem-mu-convex}
For each fixed $A\in\Rm$ we have $\mu_e(0;A)=0$, and the function
$\mu_e(\lambda;A)\ge 0$ is monotonically increasing and convex in
$\lambda\ge 0$.
\end{prop}

Proposition \ref{lem-mu-convex} allows us to define the inverse
function $\lambda_e(\mu;A)$ to $\mu_e(\lambda;A)$ that is increasing
and concave in $\mu\ge 0$ for a fixed $A$. The eigenvalue problem
(\ref{2-eig2}) may be now re-formulated as follows: given $\mu\ge 0$
find $\lambda=\lambda_e(\mu;A)$ so that the problem
\begin{equation}\label{2-eig3}
\Delta\psi-Au\cdot\nabla\psi=\mu\psi,
\end{equation}
has a solution $\psi\in E_{e,\lambda}^+$.    The variational
principle (\ref{2-c-min-mu}) for the minimal front speed now becomes
\begin{equation}\label{2-c-min-lambda}
c^*_e(A)=\inf_{\mu>0}\frac{f'(0)+\mu}{\lambda_e(\mu;A)}.
\end{equation}

\subsection{The proof of Theorem~\ref{thm-main} }

The proof of Theorem~\ref{thm-main} is based on the ideas of
\cite{NR} where they were used to obtain the asymptotics of the
pulsating front speed for cellular flows in the direction (1,0).
Here we extend them to general flows and directions $e$, and show
that they yield the conclusion of the theorem. The main ingredient
is the following Lemma.

\begin{lemma}\label{lem-lambda-mu}
There exists $\rho>0$ such that when $\eps\in[0,\tfrac 12]$ and
$\mu_0(\eps)=\rho\eps^4$, then for any $\mu\in[0, \mu_0(\eps)]$, any
flow $u$ as in Theorem \ref{thm-main}, any unit vector $e\in\Rm^2$,
and any $A\in\Rm$ we have
\begin{equation} \label{eq-lambda-mu}
(1-\eps)\frac{\sqrt\mu}{\sqrt{D_e(A)}} \le \lambda_e(\mu;A)  \le
(1+\eps)\frac{\sqrt\mu}{\sqrt{D_e(A)}}.
\end{equation}
\end{lemma}

We postpone the proof of this Lemma and prove Theorem~\ref{thm-main}
first.
Fix any $\eps>0$. Then the variational principle for $c^*_e(A)$
implies that
\begin{equation}\label{mu-l-low}
c^*_e(A)=\inf_{\mu>0}\frac{f'(0)+\mu}{\lambda_e(\mu;A)}
\ge\min\left\{\inf_{0<\mu\le\mu_0(\eps)}\frac{f'(0)+\mu}{\lambda_e(\mu;A)},
\inf_{\mu\ge\mu_0(\eps)}\frac{\mu}{\lambda_e(\mu;A)}\right\}
\end{equation}
and
\begin{equation}\label{mu-l-up}
c^*_e(A)
\le\inf_{0<\mu\le\mu_0(\eps)}\frac{f'(0)+\mu}{\lambda_e(\mu;A)}.
\end{equation}
In addition, as the function $\lambda_e(\mu,A)$ is concave and
increasing in $A$ and $\lambda_e(0;A)=0$, we have
\[
\inf_{\mu\ge\mu_0(\eps)}\frac{\mu}{\lambda_e(\mu;A)}=
\frac{\mu_0(\eps)}{\lambda_e(\mu_0(\eps);A)}\ge \sqrt{D_e(A)}
\frac{\sqrt{\mu_0(\eps)}}{1+\eps}.
\]
It follows now from (\ref{mu-l-low}) and
\eqref{eq-lambda-mu} that
\begin{equation}\label{pf-c1}
c^*_e(A)\ge \frac{\sqrt{D_e(A)}}{1+\eps}
\min\left\{\inf_{0<\mu\le\mu_0(\eps)}\frac{f'(0)+\mu}{\sqrt{\mu}},
{\sqrt{\mu_0(\eps)}}\right\}. 
\end{equation}
Using  \eqref{mu-l-up} and
\eqref{eq-lambda-mu} we also arrive at
\begin{equation}\label{pf-c2}
c^*_e(A) \le
\frac{\sqrt{D_e(A)}}{1-\eps}\inf_{0<\mu\le\mu_0(\eps)}\frac{f'(0)+\mu}{\sqrt{\mu}}.
\end{equation}

If now $f'(0)\le \tfrac 14\mu_0(\tfrac 12)$, we take
$\eps=\mu_0^{-1}(4f'(0))\le \tfrac 12$ to obtain
\[
\inf_{0<\mu\le\mu_0(\eps)}\frac{f'(0)+\mu}{\sqrt{\mu}} =
2\sqrt{f'(0)}= \sqrt{\mu_0(\eps)},
\]
and so we get from (\ref{pf-c1}) and (\ref{pf-c2}) that
\begin{equation}\label{cede}
2\sqrt{f'(0)} \frac{\sqrt{D_e(A)}}{1+\eps}\le c^*_e(A) \le
2\sqrt{f'(0)} \frac{\sqrt{D_e(A)}}{1-\eps}.
\end{equation}
Therefore, we have
\[
\left|\frac{c^*_e(A)}{2\sqrt{f'(0)D_e(A)}}-1\right|\le
\frac{\eps}{1-\eps^2}\le C(f'(0))^{1/4}
\]
for some $C>0$. This proves the second claim in Theorem
\ref{thm-main}.

On the other hand, if $f'(0)> \tfrac 14\mu_0(\tfrac 12)$, then we
set $\eps = \tfrac 12$ and obtain
\[
\sqrt{\mu_0(1/2)} \frac{\sqrt{D_e(A)}}{1+1/2}\le c^*_e(A) \le
\frac{f'(0)+\mu_0(1/2)}{\sqrt{\mu_0(1/2)}}
\frac{\sqrt{D_e(A)}}{1-1/2}.
\]
This and \eqref{cede} with $\eps\le \tfrac 12$ for $f'(0)\le \tfrac
14\mu_0(\tfrac 12)$ proves the first claim in
Theorem~\ref{thm-main}. $\Box$

\subsection{The proof of Lemma \ref{lem-lambda-mu}}

The proof of Lemma \ref{lem-lambda-mu} is similar to that of Theorem
3.1 in \cite{NR} but our result is slightly sharper.  First, we
rewrite the eigenvalue problem (\ref{2-eig3}) in terms of the
function $\zeta(x_1,x_2)=\ln\psi(x_1,x_2)$:
\begin{eqnarray}\label{eig-zeta}
&&\Delta\zeta-Au\cdot\nabla\zeta=\mu-|\nabla\zeta|^2,\\
&&\zeta(x_1+1,x_2)=\zeta(x_1,x_2) + \lambda_e(\mu;A) e_1 \nonumber\\
&&\zeta(x_1,x_2+1)=\zeta(x_1,x_2) + \lambda_e(\mu;A) e_2,\nonumber
\end{eqnarray}
where $e=(e_1,e_2)$. Without loss of generality we may assume that
$\zeta(x_1,x_2)$ has mean zero on $\Tm^2$ (which is now viewed as a
subset of $\Rm^2$):
\[
\int_{\Tm^2}\zeta(x_1,x_2)dx_1dx_2=0.
\]
As $\nabla\zeta$ is periodic both in $x_1$ and $x_2$, and
$\int_{0}^1 u_1(x_1,x_2) dx_2=0$ for each $x_1$ and $\int_{0}^1
u_2(x_1,x_2) dx_1=0$ for each $x_2$ (because $u$ is 1-periodic,
mean-zero, and incompressible), by integrating (\ref{eig-zeta}) over
$\Tm^2$ we obtain
\begin{equation}\label{mu-nabla}
\mu=\int_{\Tm^2}|\nabla\zeta|^2dx_1dx_2.
\end{equation}
The Poincar\'e inequality then implies
\begin{equation}\label{zeta-l2}
\int_{\Tm^2}|\zeta|^2dx_1dx_2\le C\mu.
\end{equation}

The function $\zeta(x_1,x_2)$ may be decomposed as
\[
\zeta(x_1,x_2)=\lambda_e(\mu;A)\left(x\cdot e
-\frac{1}{2}-\chi_e(x_1,x_2) \right) + S(x_1,x_2).
\]
Here $\chi_e(x_1,x_2)$ is the mean-zero periodic solution of
\eqref{ed-3}
and $S(x_1,x_2)$ is a mean-zero periodic correction that we would
like to show to be ``small''. Let us set
\[
\Phi(x_1,x_2)=\lambda_e(\mu;A)\left(x\cdot e -
\frac{1}{2}-\chi_e(x_1,x_2)\right).
\]
This function satisfies
\begin{eqnarray}\label{Phi-eq}
&&\Delta\Phi-Au\cdot\nabla\Phi=0\\
&&\Phi(x_1+1,x_2)=\Phi(x_1,x_2)+\lambda_e(\mu;A) e_1 \nonumber\\
&&\Phi(x_1,x_2+1)=\Phi(x_1,x_2)+\lambda_e(\mu;A) e_2.\nonumber
\end{eqnarray}
The definitions of $\Phi$ and $D_e(A)$, periodicity of $\chi_e$, and
$|e|=1$ imply that
\begin{equation}\label{phi-mu}
\|\nabla\Phi\|_2=\lambda_e(\mu;A)\sqrt{D_e(A)}.
\end{equation}

\begin{lemma}\label{lem-phi-infty}
There exists a universal constant $C>0$ such that for any flow $u$
as in Theorem \ref{thm-main}, any unit vector $e\in\Rm^2$, and any
$A\in\Rm$, we have
\begin{equation}\label{eq-phi-infty}
\|\Phi\|_{L^\infty(\Tm^2)}\le C\lambda_e(\mu;A)\sqrt{D_e(A)}
\end{equation}
and
\begin{equation}\label{eq-zeta-infty}
\zeta(x_1,x_2)\le C\sqrt{\mu} \qquad\text{for $(x_1,x_2)\in\Tm^2$}.
\end{equation}
\end{lemma}

We postpone the proof and first finish that of
Lemma~\ref{lem-lambda-mu}. It follows from (\ref{eq-phi-infty}) and
(\ref{eq-zeta-infty}) that
\begin{equation}\label{S-bd-infty}
S(x_1,x_2)\le
C\left[\lambda_e(\mu;A)\sqrt{D_e(A)}+\sqrt{\mu}\right].
\end{equation}
Note that this is only a bound from above, as is the one in
\eqref{eq-zeta-infty}. On the other hand, as
\[
\mu=\int_{\Tm^2}|\nabla\zeta|^2dx_1dx_2=\int_{\Tm^2}|\nabla(\Phi+S)|^2dx_1dx_2,
\]
the triangle inequality implies that
\begin{equation}\label{triangle}
\left|\sqrt{\mu}-\|\nabla\Phi\|_2\right|\le\|\nabla S\|_2.
\end{equation}
It follows from (\ref{eig-zeta}) and (\ref{Phi-eq}) that the
function $S$ is a mean-zero periodic solution of
\[
\Delta S-Au\cdot\nabla S=\mu-|\nabla(\Phi+S)|^2.
\]
Multiplying both sides by $S$ and integrating over $\Tm^2$ we obtain
using (\ref{S-bd-infty}),
\begin{align*}
\int_{\Tm^2}|\nabla S|^2dx_1dx_2=\int_{\Tm^2}
S|\nabla(\Phi+S)|^2dx_1dx_2 & \le
C\left[\lambda_e(\mu;A)\sqrt{D_e(A)}+\sqrt{\mu}\right] \int_{\Tm^2}
|\nabla(\Phi+S)|^2dx_1dx_2
\\ & =C\left[\lambda_e(\mu;A)\sqrt{D_e(A)}+\sqrt{\mu}\right]\mu,
\end{align*}
with $C$ from Lemma \ref{lem-phi-infty}. Using this inequality,
(\ref{triangle}), \eqref{phi-mu}, and $\sqrt{ab+c}\le\delta
a+\delta^{-1}b+\sqrt c$ with $\delta=\eps /2C$, we obtain for any
$\eps>0$
\[
\sqrt{\mu}\le \lambda_e(\mu;A)\sqrt{D_e(A)}+
C\sqrt{\lambda_e(\mu;A)\sqrt{D_e(A)}\mu+{\mu}^{3/2}}\le
\bigg(1+\frac\eps
2\bigg)\lambda_e(\mu;A)\sqrt{D_e(A)}+\frac{2C^2}{\eps}\mu+C\mu^{3/4}.
\]
We now let $\mu_0(\eps)=(\eps/4C)^4$ for $\eps\in[0, \tfrac 12]$ so
that
\[
\frac{2C^2}{\eps}\mu+C\mu^{3/4} \le \frac{3\eps}8 \sqrt\mu
\]
for all $\mu\in[0,\mu_0(\eps)]$. Then for $\mu\in[0,\mu_0(\eps)]$ we
have
\[
\sqrt{\mu}\le  (1+\eps)\lambda_e(\mu;A)\sqrt{D_e(A)}.
\]
In a similar manner, from
\[
\sqrt{\mu}\ge  \lambda_e(\mu;A)\sqrt{D_e(A)}-
C\sqrt{\lambda_e(\mu;A)\sqrt{D_e(A)}\mu+{\mu}^{3/2}}
\]
we obtain for $\mu\in[0,\mu_0(\eps)]$,
\[
(1-\eps)\lambda_e(\mu;A)\sqrt{D_e(A)}\le \sqrt{\mu}.
\]
This finishes the proof of Lemma \ref{lem-lambda-mu}. $\Box$

\subsection{The proof of Lemma \ref{lem-phi-infty}}

Both statements of this lemma are an immediate consequence of the
following two propositions.
\begin{prop}\label{prop-Linfty}
There exists a universal constant $C>0$ such that if
$\alpha\in\Rm^2$ and the function $q$ on $\Rm^2$ satisfies
\begin{eqnarray*}
&& \Delta q-A  u \cdot \nabla q+  |\nabla q|^2 \ge 0,\\
&& q(x_1+1,x_2) = q(x_1,x_2) + \alpha_1,\\
&& q(x_1, x_2+1) =q(x_1,x_2)+ \alpha_2,\\
&& \int_{\Tm^2}q(x_1,x_2)\,dx_1dx_2=0,
\end{eqnarray*}
then
\begin{equation}\label{L-infty-q}
q(x_1,x_2)\le C \big(\|\nabla
q\|_{L^2({\Tm^2})}+|\alpha_1|+|\alpha_2| \big)
\end{equation}
for all $(x_1,x_2)\in\Tm^2$.
\end{prop}

\noindent{\it Remark.} When $\alpha_2=0$ then this is Proposition 4
in \cite{NR}. \medskip

\noindent{\bf Proof.} The Poincar\' e inequality says that
\begin{equation}\label{Poinc}
\|q\|_{L^2({\Tm^2})}\le \bigg( \frac C2 -1 \bigg) \|\nabla
q\|_{L^2({\Tm^2})}
\end{equation}
for some $C>2$. We will show that \eqref{L-infty-q} holds with this
$C$.

Let
\[
M:=\max_{(x_1,x_2)\in{\Tm^2}}q(x_1,x_2) = q(\vx_0)
\]
for some $\vx_0\in\Tm^2$ and consider the case
\[
M\ge 2\|q\|_{L^2({\Tm^2})}+|\alpha_1|+|\alpha_2|
\]
(otherwise we are done by \eqref{Poinc}).

Assume first that
\begin{equation}\label{3-exists-large-q}
\hbox{for any }x_1\in[0,1]\hbox{ there exists }s(x_1)\in[0,1]\hbox{
such that } q(x_1,s(x_1))\ge M-|\alpha_1|-|\alpha_2|.
\end{equation}
Observe that if we define the set ${\cal A}\subset[0,1]$ as
\begin{equation}\label{3-define-A-q}
{\cal A}=\left\{x_1\in[0,1]:~\exists r(x_1)\in[0,1] \hbox{ such that
}|q(x_1,r(x_1))|\le 2\|q\|_{L^2({\Tm^2})}\right\}
\end{equation}
then
\begin{equation}\label{meas-good-q}
\left|{\cal A}\right|\ge \frac{1}{2}.
\end{equation}
Now, (\ref{3-exists-large-q})--(\ref{meas-good-q}) imply that
\[
\int_0^1 \left|\pdr{q}{x_2}\right|^2dx_2\ge  \bigg[ \int_0^1
\pdr{q}{x_2}\,dx_2 \bigg]^2 \ge \left(M- |\alpha_1|-|\alpha_2| -
2\|q\|_{L^2({\Tm^2})}\right)^2
\]
for all $x_1\in\cal A$. We then obtain
\[
\int_{\Tm^2}|\nabla q|^2{dx_1dx_2}\ge \int_{\cal
A}\left|\pdr{q}{x_2}\right|^2{dx_1dx_2}\ge \frac{1}{2}\left(M
-|\alpha_1|-|\alpha_2| - 2\|q\|_{L^2({\Tm^2})}\right)^2.
\]
It follows from the above and \eqref{Poinc} that
\[
\max_{(x_1,x_2)\in{\Tm^2}} q(x_1,x_2) =M \le 2\left( \|\nabla
q\|_{L^2({\Tm^2})} + \|q\|_{L^2({\Tm^2})}
+|\alpha_1|+|\alpha_2|\right)\le C\left(\|\nabla
q\|_{L^2({\Tm^2})}+|\alpha_1|+|\alpha_2|\right).
\]

If \eqref{3-exists-large-q} does not hold, but
\begin{equation}\label{3-exists-large-q2}
\hbox{for any }x_2\in[0,1]\hbox{ there exists }s(x_2)\in[0,1]\hbox{
such that } q(s(x_2),x_2)\ge M-|\alpha_1|-|\alpha_2|,
\end{equation}
then an identical argument applies. Assume therefore that neither
\eqref{3-exists-large-q} nor \eqref{3-exists-large-q2} hold. That
means that there are $\xi_1,\xi_2\in[0,1]$ such that for all
$\xi\in[0,1]$,
\[
q(\xi_1,\xi),q(\xi,\xi_2)< M-|\alpha_1|-|\alpha_2|.
\]
If now $S_j$, $j=1,2,3,4$, are the four unit squares with sides
parallel to the axes, and having $(\xi_1,\xi_2)$ as a vertex, then
for each $j$ we must have
\[
\max_{(x_1,x_2)\in\partial S_j} q(x_1,x_2) < M-|\alpha_1|-|\alpha_2|
+ |\alpha_1|+|\alpha_2|=M
\]
by the assumptions on $q$. But one of these squares contains $\vx_0$
with $q(\vx_0)=M$, which is a contradiction because $q$ satisfies
the maximum principle. Hence \eqref{3-exists-large-q} or
\eqref{3-exists-large-q2} hold and we are done. $\Box$

\begin{prop}\label{lem-lambda-mu-sqrt}
For all $A$ we have
\begin{equation}\label{eq-lem-lambda-mu-sqrt}
\lambda_e(\mu;A)\le\sqrt{\mu}.
\end{equation}
\end{prop}

\noindent{\bf Proof.} This is an immediate consequence of
\eqref{eig-zeta}, (\ref{mu-nabla}), and $|e|=1$:
\[
\mu=\int_{\Tm^2}|\zeta_{x_1}|^2 + |\zeta_{x_2}|^2 dx_1dx_2\ge
\left(\int_{\Tm^2}\zeta_{x_1}dx_1dx_2\right)^2  +
\left(\int_{\Tm^2}\zeta_{x_2}dx_1dx_2\right)^2 =\lambda_e(\mu;A)^2.
\quad\Box
\]

\noindent{\bf Proof of Lemma~\ref{lem-phi-infty}.} This holds with
$C$ being $(1+\sqrt 2)$-times the constant from
Proposition~\ref{prop-Linfty}. The upper bound on $\Phi(x_1,x_2)$
follows directly from (\ref{Phi-eq}), \eqref{phi-mu}, and $D_e(A)\ge
1$. The lower bound on $\Phi$ follows after applying
Proposition~\ref{prop-Linfty} to $-\Phi$. The upper bound
(\ref{eq-zeta-infty}) for $\zeta(x_1,x_2)$ follows from
Propositions~\ref{prop-Linfty} and \ref{lem-lambda-mu-sqrt}, and
(\ref{mu-nabla}). $\Box$

\section{An extension to higher dimensions and examples}\label{sec:gen}

\subsection{The main result in higher dimensions}

As mentioned before, Proposition \ref{thm-main0} holds in any
dimension. Moreover, the following partial analog of Corollary
\ref{cor-main} holds.

\begin{thm} \label{T.4.1}
Let $u(x)$ be a 1-periodic incompressible Lipschitz flow on $\Rm^n$,
let $f$ be a KPP nonlinearity, and $e\in\Rm^n$ a unit vector. If
\eqref{eq-thm-eigen-fn} has no $H^1(\Tm^n)$-solutions, then
\begin{equation}\label{4.1}
\lim_{A\to+\infty}c^*_e(Au,f)=\lim_{A\to+\infty}D_e(Au,f)=+\infty.
\end{equation}
\end{thm}

\noindent{\bf Proof.} Assume that there is a sequence
$A_k\to+\infty$ such that
\begin{equation}\label{4.3}
M:=\sup_nc^*_e(A_k)<+\infty.
\end{equation}
Let $\varphi_k$ be a solution of \eqref{2-eig1} on $\Tm^n$ with
$A=A_k$ and $\lambda=\lambda_k$ where $\lambda_k$ is such that
\begin{equation}\label{4.4}
c^*_e(A_k)=\frac{f'(0)+\lambda_k^2+\kappa_e(\lambda_k;A_k)}{\lambda_k}
\end{equation}
with $\kappa_e$ defined in \eqref{2-eig1}. Moreover, we choose
$\varphi_k$ so that $\omega_k=\ln\varphi_k$ is mean-zero on $\Tm^n$.
Denote $\kappa_k=\kappa_e(\lambda_k;A_k)$ so that
\begin{equation}\label{4.5}
\Delta\omega_k+|\nabla\omega_k|^2 - A_ku\cdot\nabla\omega_k
-2\lambda_k e\cdot\nabla \omega_k + A_k\lambda_k u\cdot e =\kappa_k.
\end{equation}
Integrate this over $\Tm^n$ and use that $u$ is incompressible and
mean-zero to obtain
\[
\kappa_k= \int_{\Tm^n} |\nabla\omega_k|^2\,dx\ge 0.
\]
This, \eqref{4.3}, and \eqref{4.4} mean that the $\lambda_k$ are
bounded away from 0 and $\infty$ (namely, $|2\lambda_k-M|\le
\sqrt{M^2-4f'(0)}\,$). So after passing to a subsequence we can
assume $\lambda_k\to\lambda_0\in(0,\infty)$. This and \eqref{4.3}
give
\[
\sup_k \|\nabla\omega_k\|^2_{L^2(\Tm^n)} = \sup_k \kappa_k<+\infty.
\]
By \eqref{Poinc}, the $\omega_k$ are bounded in $H^1(\Tm^n)$, and so
converge to some $\omega_0\in H^1(\Tm^n)$ strongly in $L^2(\Tm^n)$
and weakly in $H^1(\Tm^n)$. Then also
$\Delta\omega_k\to\Delta\omega_0$ in the sense of distributions and
$|\nabla\omega_k|^2$ are bounded in $L^1(\Tm^n)$. Divide \eqref{4.5}
by $A_k$ and pass to the limit so that in the sense of
distributions,
\[
- u\cdot\nabla\omega_0  + \lambda_0 u\cdot e = 0.
\]
Since $\omega_0\in H^1(\Tm^n)$, we now see that
$\phi_e=\omega_0/\lambda_0\in H^1(\Tm^n)$ solves
\eqref{eq-thm-eigen-fn}, yielding a contradiction. $\Box$

\subsection{Examples}

We conclude with several examples of flows in two and three
dimensions which speed up reaction-diffusion fronts. We will
consider the case $e=e_1=(1,0)$ (resp. $e=e_1=(1,0,0)$ in three
dimensions) so that \eqref{eq-thm-eigen-fn} becomes
\begin{equation} \label{4.6}
u\cdot\nabla\phi=u_1,
\end{equation}
but our analysis easily extends to other directions 
of front propagation.

\begin{example} \label{E.4.1a}{\bf Shear flows.}
These are flows of the form $u(x_1,x_2)=(\alpha(x_2),0)$ with
mean-zero $\alpha$. In this case \eqref{2-eig1} with $e=e_1$ becomes
\[
\Delta \varphi-(A\alpha(x_2)+2\lambda)\varphi_{x_1} +
{\lambda}A\alpha(x_2)\varphi=\kappa(\lambda;A)\varphi,\qquad
\hbox{$\varphi>0$}.
\]
The unique solution to this equation is of the form
$\varphi(x_1,x_2)=\varphi(x_2)$ and satisfies
\begin{equation} \label{4.7a}
\varphi_{x_2x_2} +
{\lambda}A\alpha(x_2)\varphi=\kappa(\lambda;A)\varphi,\qquad
\hbox{$\varphi>0$}.
\end{equation}
Obviously then $\kappa(\lambda;A)=\kappa(\lambda A;1)$ and so we
have
\[
\lim_{A\to\infty} \frac{c^*_e(A)}A = \lim_{A\to\infty}
\inf_{\lambda>0}\frac{f'(0)+\lambda^2+\kappa(\lambda A;1)}{\lambda
A} = \lim_{A\to\infty}
\inf_{\lambda>0}\frac{f'(0)+\frac{\lambda^2}{A^2}+\kappa(\lambda;1)}{\lambda}
= \inf_{\lambda>0}\frac{f'(0)+\kappa(\lambda;1)}{\lambda}.
\]
Notice that this shows that $c^*_e(A)/A$ is decreasing to a positive
limit. This has been established in \cite{B} and an alternative
(variational) characterization of the limit has been provided in
\cite{Heinze}. Multiplying \eqref{4.7a} by $\varphi$ and integrating
over $\Tm^2$ one obtains $\kappa(\lambda;1)\le
\lambda\|\alpha_+\|_\infty$ with $\alpha_+$ the positive part of
$\alpha$. Since the operator on the RHS of \eqref{4.7a} is
self-adjoint, it is easy to see that in fact
$\lim_{\lambda\to\infty} \kappa(\lambda;1)/\lambda =
\|\alpha_+\|_\infty$. From this it follows that for each $f'(0)$,
\[
\lim_{A\to\infty} \frac{c^*_e(A)}A \le \|\alpha_+\|_\infty =
\lim_{f'(0)\to\infty} \lim_{A\to\infty} \frac{c^*_e(A)}A
\]
On the other hand, \eqref{ed-3} becomes
\[
-\chi_{x_2x_2}=A\alpha(x_2)
\]
with $\chi(x_1,x_2)=\chi(x_2)$, and it follows that
$D_e(A)=1+A^2\gamma^2$ with
$\gamma=\|\nabla_{x_2}(-\Delta_{x_2})^{-1}\alpha\|_2$. We thus have
\[
\lim_{A\to\infty} \frac{c^*_e(A)}{\sqrt{D_e(A)}} \le
\frac{\|\alpha_+\|_\infty}\gamma
\]
for any $f'(0)$ and $\alpha\not\equiv 0$. This shows that the lower
bound in \eqref{intro-D-c} is optimal up to a constant for large
$f'(0)$. Note that this example extends without change to any
dimension.
\end{example}

If a solution $\phi\in H^1(\Tm^2)$ of \eqref{4.6} exists, then the
function $\psi(x_1,x_2)= \phi(\{x_1\},x_2)-x_1$ belongs to
$H^1(2\Tm\times\Tm,2\Tm)$ (here $2\Tm=[0,2]$ with 0 and 2
identified, and $\{x\}$ is the fractional part of $x$; in three
dimensions  $\psi\in H^1(2\Tm\times\Tm^2,2\Tm)$). Moreover, $\psi$
satisfies
\begin{equation} \label{4.7}
u\cdot\nabla\psi=0 \qquad\text{and}\qquad
\psi(x_1+1,x_2)=\psi(x_1,x_2)-1.
\end{equation}
In the following examples it will be a bit easier to work with
$\psi$ than with $\phi$.

\begin{example} \label{E.4.2} {\bf Percolating flows.}
Let $u$ be such that for some $x_2$ the streamline $\calS$ of the
flow (i.e., solution of the ODE $X'=u(X)$) starting at $(0,x_2)$
reaches $(1,x_2)$ in finite time. This means that $\calS$ is a
periodic curve passing through  $(n,x_2)$ for each $n\in\Zm$. Such
$u$ is called a percolating flow (in the horizontal direction). Note
that any non-zero shear flow is percolating. Since $u\neq 0$ on
$\calS$, it is easy to show that if $\psi$ from \eqref{4.7} exists,
it must be continuous and constant on $\calS$ (see \cite[Lemma
5.2]{ZlaMix}). This, however, contradicts the second condition in
\eqref{4.7}. Hence such a $\psi$ does not exist and Theorem
\ref{thm-main} gives
\begin{equation}\label{example-infinite}
\lim_{A\to+\infty}c_e^*(A)=+\infty
\end{equation}
for $e=e_1$.
This has been proved under some additional ``non-degeneracy''
conditions on $u$ in \cite{CKOR}.
\end{example}

\begin{example} \label{E.4.3} {\bf Cellular flows.}
These are flows with a periodic array of cells, each streamline of
the flow being either a closed curve contained inside a cell or a
part of the boundary of a cell. A typical example is the flow
\begin{equation} \label{4.8}
u(x_1,x_2)=\nabla^\perp
H(x_1,x_2)=(-H_{x_2}(x_1,x_2),H_{x_1}(x_1,x_2))
\end{equation}
with $H(x_1,x_2)=\sin 2\pi x_1 \sin 2\pi x_2$ the stream function of
$u$. The streamlines of this flow are depicted in
Figure~\ref{fig-cell}.

\begin{figure}[ht!]
 \centerline{\epsfxsize=0.6\hsize \epsfbox{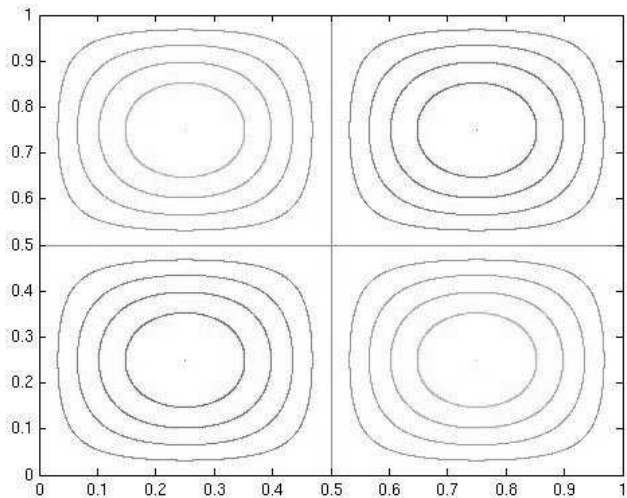}} 
\caption{A 2D cellular flow.}
 \label{fig-cell}
\end{figure}

Constancy and continuity of $\psi$ from \eqref{4.7} on each
non-trivial streamline implies that $\psi$ has to be constant on the
boundary of each cell (namely, it equals there the limit of the values
on the streamlines approaching the boundary). But then $\psi$ can
belong to $H^1$ only if it is constant on the whole ``skeleton'' of
separatrices separating the cells, which again contradicts the second
condition in \eqref{4.7}. Hence \eqref{example-infinite} holds for
$e=e_1$. We note that front speed-up by (certain generic) cellular
flows has first been proved in \cite{KR}, with precise asymptotics
established in \cite{NR}.
\end{example}

\begin{example} \label{E.4.4} {\bf Checkerboard flows.}
Consider the cellular flows from the previous example with the flow
vanishing in every other cell, thus forming the checkerboard pattern
 depicted in Figure~\ref{fig-checker}.

\begin{figure}[ht!]
 \centerline{\epsfxsize=0.6\hsize \epsfbox{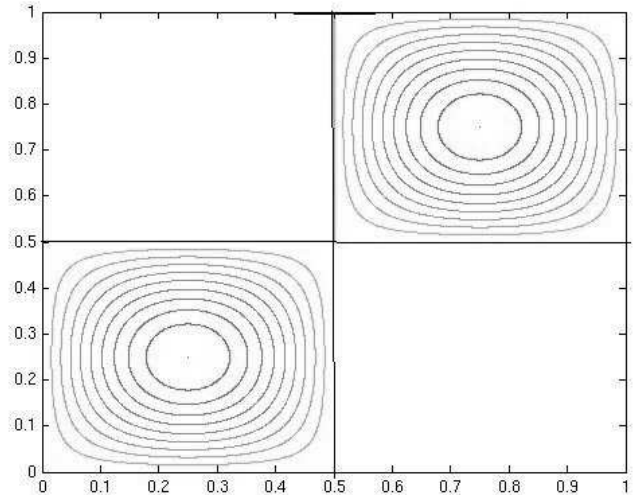}} 
 \caption{A checkerboard cellular flow.}
 \label{fig-checker}
\end{figure}

The above $u$ is not Lipschitz in this case, so let us remedy this
problem by taking, for instance, $H(x_1,x_2)=(\sin 2\pi x_1 \sin
2\pi x_2)^\alpha$ with $\alpha\ge 2$ in the cells where $u$ does not
vanish. Now $u$ vanishes on the boundaries of all cells, but if
$\psi$ from \eqref{4.7} exists, the requirement $\psi\in H^1$ still
ensures that $\psi$ is constant on the boundaries of those cells in
which the flow does not vanish (and it is continuous at these
boundaries from inside of these cells). The values of these
constants for two cells that touch by a corner must be the same.
This follows from the fact that if an $H^1$ function has values $a$
and $b$ on two curves connecting a point $x\in\Rm^2$ to a circle
$\{y:~|y-x|=\eps\}$
for some $\eps>0$ (and it is continuous at these curves,
at least from one side), then $a=b$ (see, e.g. \cite[Lemma
5.2]{ZlaMix} for this simple argument). We again have contradiction
with the second condition in \eqref{4.7}, and
\eqref{example-infinite} for $e=e_1$ follows. 
\end{example}

Effective diffusivity enhancement for the
flow in Figure~\ref{fig-checker} has been proved in \cite{FP}. 
Notice, however,
that our argument above can handle more general checkerboard-type flows.
Consider, for instance, the case of such flows but with the contact
between cells in which the flow does not vanish formed only by two
touching cusps rather than two right angles. In this case the angle
of contact between the cells can be equal to $\pi$ but we still have
traveling front speed enhancement by the flow in the sense of
\eqref{example-infinite}.

\begin{example} \label{E.4.5} {\bf Flows with gaps.}
Consider the
cellular flow from Example \ref{E.4.3} but with a vertical ``gap''
of width $\delta>0$ (in which $u=0$) inserted in place of each
vertical line $\{n\}\times\Rm$, $n\in\Zm$, such as shown in
Figure~\ref{fig-gaps}.

\begin{figure}[ht!]
 \centerline{\epsfxsize=0.6\hsize \epsfbox{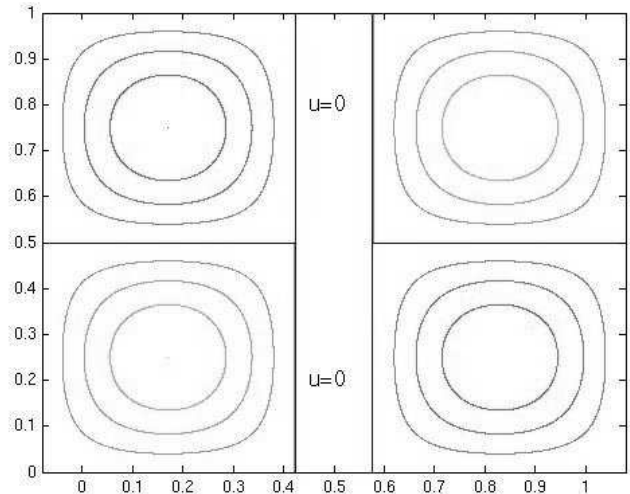}} 
 \caption{A cellular flow with gaps.}
 \label{fig-gaps}
\end{figure}

This can be achieved by
letting the stream function be, for instance,
\[
H(x_1,x_2)=
\begin{cases}
(\sin \frac{2\pi}{1-\delta} x_1 \sin 2\pi x_2)^\alpha & \quad
x_1\text{ mod } 1 \in[0,1-\delta)
\\ 0 & \quad x_1 \text{ mod } 1 \in[1-\delta,1) \end{cases}
\]
with $\alpha\ge 2$. More generally, we can assume that $u$ has a
periodic array of vertical channels in which the flow only moves
``along'' each channel. Namely, we assume that there is a $C^1$ map
$\gamma:[0,1]\times\Rm\to [0,1]\times\Rm$ which is one-to-one and
satisfies $\gamma(a,b+1)=\gamma(a,b)+(0,1)$ and $u(\gamma(a,b))=
c(a,b) \tfrac{\partial}{\partial y}\gamma(a,b)$ for all $(a,b)\in
[0,1]\times\Rm$ and some $c(a,b)\in \Rm$. This means that
$\gamma(\{a\}\times\Rm)$ is a ``vertical'' streamline of $u$ for
each $a$ (or a union of streamlines if $u$ vanishes somewhere on
this curve). Let us also assume that the curve
$\gamma(\{0\}\times\Rm)$ lies to the ``left'' of
$\gamma(\{1\}\times\Rm)$. We now let $\omega:[0,1]\to[0,1]$ be a
$C^1$ function with $\omega(0)=1$ and $\omega(1)=0$. We define
$\psi(\gamma(a,b))=\omega(a)$ and set $\psi$ equal 1 and 0 on the
``left'' and ``right'' components of $[0,1]\times\Rm
\smallsetminus\gamma([0,1]\times\Rm)$, respectively. We have that
$\psi$ is constant on each $\gamma(\{a\}\times\Rm)$ and so
$u\cdot\nabla\psi=0$. Since the periodicity of $\gamma$ ensures
$\psi(x_1,0)=\psi(x_1,1)$, we can extend $\psi$ to an
$H^1(2\Tm\times\Tm,2\Tm)$-function by letting
$\psi(x_1+1,x_2)=\psi(x_1,x_2)-1$. Hence \eqref{4.7} is satisfied
and \eqref{example-infinite} fails for $e=e_1$.
\end{example}

Note that the flows in the last example are percolating in the
vertical direction if $u$ does not vanish on each of the curves
$\gamma(\{a\}\times\Rm)$. Hence by the remark after
Proposition~\ref{thm-main0} and Example~\ref{E.4.2}, speed-up of
fronts in the sense of \eqref{example-infinite} happens precisely
when $e\neq e_1$. On the other hand, \eqref{example-infinite} is 
valid for any $e$ in the case of the flows from Examples \ref{E.4.3} and \ref{E.4.4}.

Our final example deals with speed-up of fronts by three-dimensional
cellular flows. To the best of our knowledge this is the first
example of this kind.

\begin{example} \label{E.4.6} {\bf 3D cellular flows.} We consider here
flows that have a cellular structure and are truly
three-dimensional, with all three components of the velocity
depending on all three coordinates. Such incompressible flows have
been constructed in \cite{Bisshopp}. They have the form
\begin{equation} \label{4.9}
u(x_1,x_2,x_3) =  \big( \Phi_{x_1}(x_1,x_2) W'(x_3),
\Phi_{x_2}(x_1,x_2) W'(x_3), k\Phi(x_1,x_2)W(x_3) \big)
\end{equation}
with $\Delta\Phi\equiv -k\Phi$. We will concentrate here on the
simplest example of a flow with cubic cells. This flow is given by
\begin{equation} \label{4.10}
\Phi(x_1,x_2) = \cos x_1\cos x_2, \qquad W(x_3) = \sin x_3
\end{equation}
and $k=2$, with $W$ possibly any other $2\pi$-periodic function
vanishing at $0$.

The cube $\calC=[0,2\pi]^3=(2\pi\Tm)^3$ is a cell of periodicity for
the flow \eqref{4.10}. Each of the eight dyadic sub-cubes of $\calC$
of side-length $\pi$ is invariant under the flow and the flow in
each of them has the same streamlines but moves along them in
opposite directions in two neighboring sub-cubes. For the moment we
restrict our attention to only one of them, $[0,\pi]^3$. The
streamlines of the flow in the part of this sub-cube given by
$x_1+x_2\ge \pi$ are depicted in Figure~\ref{fig-cell3d} (and they are
symmetric across the plane $x_1+x_2=\pi$).

\begin{figure}[ht!]
\centerline{\epsfxsize=0.6\hsize \epsfbox{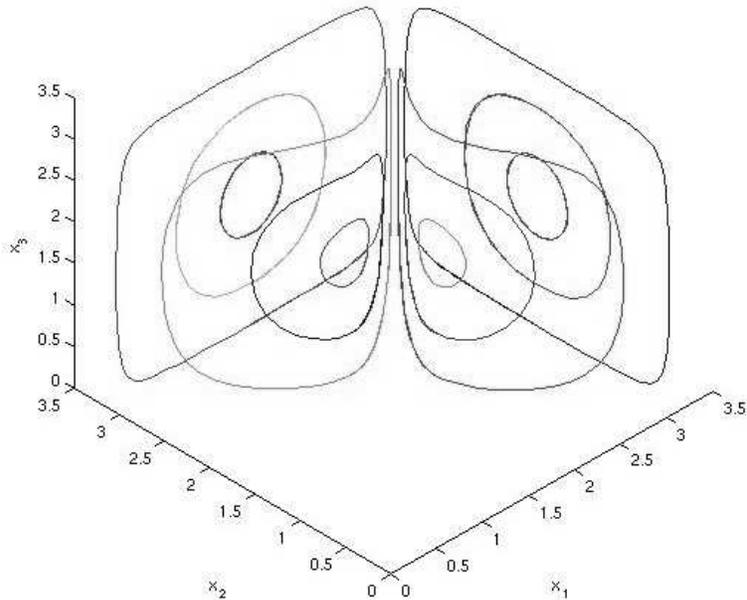} }
 \caption{A 3D cellular flow.}
 \label{fig-cell3d}
\end{figure}

Since
\[
\frac{u_1(x_1,x_2,x_3)}{u_2(x_1,x_2,x_3)} =
\frac{\Phi_{x_1}(x_1,x_2)}{\Phi_{x_2}(x_1,x_2)} = \frac{\sin x_1
\cos x_2}{\cos x_1\sin x_2}
\]
is independent of $x_3$, the projection onto the  $(x_1,x_2)$--plane
of any streamline of the flow stays on a curve satisfying
$X'=\nabla\Phi(X)$. Figure \ref{fig-2d-3d} shows the phase portrait
of this ODE in $[0,\pi]^2$.

\begin{figure}[ht!]
 \centerline{\epsfxsize=0.3\hsize \epsfbox{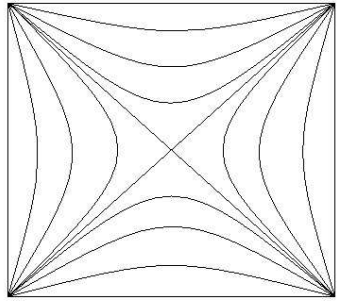}} 
 \caption{A 2D projection of the 3D cellular flow.}
 \label{fig-2d-3d}
\end{figure}

With the exception of the diagonals
$D=\{(x_1,x_2)\,:\, |x_1-\pi/2|=|x_2-\pi/2|\}$, each of the pictured
curves connects a minimum and a maximum of $\Phi$ (we will call
these $\Phi$-curves). Any streamline of $u$ in $[0,\pi]^3$, with the
exception of those lying in
$(D\times[0,\pi])\cup([0,\pi]^2\times\{0,\pi\})$ (these form a set
of measure zero and will be disregarded from now on), is a closed
orbit, whose projection onto the $(x_1,x_2)$--plane is a
portion of a $\Phi$-curve contained in $[0,\pi]^2\setminus D$ (in
particular, it does not contain the corners of $[0,\pi]^2$). Some of
these are stationary, namely those at $(x_1,x_2,x_3)$ with
$x_3=\pi/2$ and either $x_1=\pi/2$ or $x_2=\pi/2$.

For each $(x_1,x_2)\in[0,\pi]^2\setminus D$, let
$J(x_1,x_2)=(x_1,\pi-x_2)$ if $|x_1-\pi/2|>|x_2-\pi/2|$ and
$J(x_1,x_2)=(\pi-x_1,x_2)$ if $|x_1-\pi/2|<|x_2-\pi/2|$, so that
$(x_1,x_2)$ and $J(x_1,x_2)$ lie on the same $\Phi$-curve. Extend
$J$ to the whole cube $[0,\pi]^3$ by letting
$J(x_1,x_2,x_3)=(J(x_1,x_2),x_3)$ whenever
$|x_1-\pi/2|\neq|x_2-\pi/2|$. Notice that $J$ is defined for almost
all points in the cube (and we will from now on disregard the rest)
and $J^2={\rm Id}$. Each of the streamlines of $u$ is symmetric
either across the plane $x_1=\pi/2$ or across $x_2=\pi/2$, and hence
the points $x$ and $J(x)$ lie on the same streamline of $u$.

Let us now assume $\psi\in H^1(2\calC,4\pi\Tm)$ with
$2\calC=4\pi\Tm\times (2\pi\Tm)^2$ satisfies
\begin{equation} \label{4.11}
u\cdot\nabla\psi=0 \qquad\text{and}\qquad
\psi(x_1+2\pi,x_2,x_3)=\psi(x_1,x_2,x_3)-2\pi.
\end{equation}
(which is the analogue of \eqref{4.7}) and assume without loss of
generality that $\psi$ is real. Extend $J$ to $2\calC$ by
$J(x_1+p\pi,x_2+q\pi,x_3+r\pi)=J(x_1,x_2,x_3)+\pi(p,q,r)$ for any
$p\in\{0,1,2,3\}$ and $q,r\in\{0,1\}$ so that one still has that $x$
and $J(x)$ lie on the same streamline of $u$. The condition
$u\cdot\nabla\psi=0$ implies that $\psi$ is constant on almost all
streamlines of $u$. In particular, for almost all $x$ we have
$\psi(x)=\psi(J(x))$. At the same time we know that the restriction
of $\psi$ to almost every plane $x_3=C$ is an $H^1(4\pi\Tm\times
2\pi\Tm)$ function. This means that for almost every plane
$x_3=C\in(0,\pi)$  the restriction of $\psi$ to it (which we
call $\bar\psi$) belongs to $H^1(4\pi\Tm\times 2\pi\Tm)$ and satisfies
$\bar\psi(x)=\bar\psi(J(x))$ for almost all $x=(x_1,x_2)\in 4\pi\Tm\times
2\pi\Tm$.

Next we choose $r_0$ such that if  $B_j= B((j\pi,0),r_0)$, then
$\int_{B_j}|\nabla\bar\psi(x)|^2 dx < \pi/32$ for $j=0,1,2$. It is then
easy to show that for each $j=0,1,2$ the set $R_j$ of all $r\in
(0,r_0)$ such that
\[
{\rm ess} \sup_\theta\bar \psi(j\pi+r\cos\theta,r\sin\theta) - {\rm ess}
\inf_\theta \bar\psi(j\pi+r\cos\theta,r\sin\theta) < \frac \pi 2
\]
satisfies $|R_j|>3r_0/4$. This is because
\[
\frac \pi{32}> \int_{B_j} |\nabla\bar\psi|^2 dx \ge
\int_{(0,r_0)\setminus R_j} \int_0^{2\pi} r \bigg| \frac 1r
\frac{\partial\bar\psi}{\partial\theta} \bigg|^2 d\theta dr \ge
\int_{(0,r_0)\setminus R_j} \frac 1{2\pi r} \bigg( \int_0^{2\pi}
\bigg| \frac{\partial\bar\psi}{\partial\theta} \bigg| d\theta \bigg)^2
\ge \frac {\pi|(0,r_0)\setminus R_j|}{8 r_0}.
\]
Let $R_0=\bigcap_{j=0}^2 R_j$ so that $|R_0|> r_0/4$, and let $r\in
R_0$. Then the values $\psi(x)$ for almost all points $x$ on the
circle $C_j(r)=\{(j\pi+r\cos\theta,r\sin\theta) \,|\,
\theta\in[0,2\pi)\}$ lie within an interval $|I_j(r)|\le \pi/2$.
Since
$\bar\psi(r\cos\theta,r\sin\theta)=\bar\psi(\pi-r\cos\theta,r\sin\theta)$
and $\bar\psi(\pi
+r\cos\theta,r\sin\theta)=\bar\psi(2\pi-r\cos\theta,r\sin\theta)$ for
almost all $(r,\theta)\in R_0\times(-\pi/4,\pi/4)$ (because
$\bar\psi(x)=\bar\psi(J(x))$), we have that for almost all $r\in R_0$, the
values $\bar\psi(x)$ for almost all $x\in \bigcup_{j=0}^2 C_j(r)$ lie
within the interval $I_0(r)=\bigcup_{j=0}^2 I_j(r)$ with
$|I_0(r)|\le 3\pi/2$. But this contradicts the second condition in
\eqref{4.11} because $C_2(r)=C_0(r)+(2\pi,0)$. Therefore there is no 
$\psi\in H^1(2\calC,4\pi\Tm)$ which satisfies \eqref{4.11}, and hence 
\eqref{example-infinite} holds for $e=e_1$.
\end{example}

We note that our analysis can also be performed on other 3D flows
from \cite{Bisshopp}, for example, on the flow given by
\[
\Phi(x_1,x_2) = 2\cos\sqrt 3x_1\cos x_2 + \cos \bigg( 2x_2-\frac \pi
2\bigg), \qquad W(x_3)=\sin x_3,
\]
whose cells form a hexagonal 3D honeycomb lattice. Using the fact that 
$\Phi(x_1,x_2) = -\Phi(R(x_1,x_2))$, with $R$ the reflection accross 
any of the lines $x_2=kx_1-\tfrac\pi 2$, $k=\pm\sqrt{3},0$, one 
can show as above that the streamlines of the flow are symmetric accross 
the planes given by these three equations, and again conclude 
\eqref{example-infinite}.

\end{document}